\documentclass[12pt]{article}
\usepackage{latexsym}
\usepackage{amsfonts}
\def\Bbb{\mathbb}
\def\bcp{\mathbb C \mathbb P}
\def\RR{\mathbb R}

\def\ZZ{\mathbb Z}
\def\ro{\stackrel{\circ}{r}}
\def\eea{\end{eqnarray*}}
\def\mv{\mbox{Vol}}

\newtheorem{thm}{Theorem}[section]
\newtheorem{prop}[thm]{Proposition}
\newtheorem{cor}[thm]{Corollary}

\newtheorem{conj}[thm]{Conjecture}
\newtheorem{defn}[thm]{Definition}
\newenvironment{proof}{\medskip \noindent
{\bf Proof.}}{\hfill \rule{.5em}{1em}
\\}
\newenvironment{xpl}{\mbox{ }\\ {\bf  Example}\mbox{ }}{
\hfill $\diamondsuit$\mbox{}\bigskip}

\title{Curvature and Smooth Topology \\
in Dimension Four}

\author{Claude LeBrun\thanks{Supported 
in part by  NSF grant DMS-9802722.} 
\\SUNY Stony Brook}
\date{}
\begin{document}
\maketitle

  \begin{abstract}	
 Seiberg-Witten theory leads to a delicate interplay 
 between  Riemannian geometry  and smooth topology in
 dimension four. In particular, 
 the  scalar curvature of any metric must satisfy certain non-trivial
 estimates if the manifold in question has a non-trivial
 Seiberg-Witten invariant. However, it has recently been 
 discovered \cite{lebweyl,lric} that similar statements 
 also apply to other parts of the curvature tensor. 
 This article  presents the most salient aspects 
 of these curvature estimates in a self-contained manner, 
 and shows how they can be applied
 to the theory of Einstein manifolds. 
 We then probe the issue of whether the
 known estimates are optimal by  relating this question to a 
 certain conjecture in K\"ahler geometry. 
 	\end{abstract}
  
\vfill
\pagebreak

\section{Four-Dimensional Geometry}
\label{prim}

Let $M$ be a smooth compact oriented  4-manifold. 
If $g$ is any Riemannian metric on $M$, the 
middle-dimensional Hodge star operator 
$$\star: \Lambda^2 \to \Lambda^2$$
 has eigenvalues $\pm 1$. 
This gives rise to a natural decomposition 
\begin{equation} 
\Lambda^2 = \Lambda^+ \oplus \Lambda^- 
\label{deco} 
\end{equation}
of the 
rank-6 bundle of 2-forms 
into two rank-3 bundles,
where 
$$\psi \in \Lambda^{\pm}~~\Longleftrightarrow~~
\star\psi = \pm \psi.$$ 
Sections of $\Lambda^+$ are called {\em self-dual}
2-forms, while sections of 
$\Lambda^-$ are called {\em anti-self-dual}
2-forms. The middle-dimensional  Hodge star operator 
is unchanged if $g$ is multiplied by a smooth positive 
function, so the 
decomposition 
(\ref{deco}) really only depends on the conformal class
$\gamma =[g]$ rather than on the Riemannian metric itself.

The decomposition (\ref{deco}) has 
important ramifications for Riemannian geometry. 
In particular, since the 
curvature tensor $\mathcal R$ may be thought of as a linear map
$\Lambda^{2}\to \Lambda^{2}$, there is an induced  
 decomposition \cite{st}
$$
{\mathcal R}=
\left(
\mbox{
\begin{tabular}{c|c}
&\\
$W_++\frac{s}{12}$&$\ro$\\ &\\
\cline{1-2}&\\
$\ro$ & $W_{-}+\frac{s}{12}$\\&\\
\end{tabular}
} \right) 
$$
into simpler curvature tensors. 
Here the
{\em self-dual} and {\em anti-self-dual Weyl curvatures}  
$W_\pm$ are the trace-free pieces of the appropriate blocks.
The scalar curvature  $s$ is understood to act by scalar multiplication,
and $\ro$ can be identified with the trace-free part 
$r-\frac{s}{4}g$  of the Ricci curvature.

Now the $L^2$-norm of each of these curvatures
are all scale invariant, so one may define  sensible 
diffeomorphism invariants
of a 4-manifold, such as  
$$	(\inf\|s\|)(M)  =  \inf_{g} \left( \int_{M} 
	s^{2}_{g}d\mu_{g}\right)^{1/2} \\
$$
by considering the infima of the  $L^2$ norms
of these curvatures 
over all
metrics $g$ on $M$. 
Notice that these invariants  might {\em a priori} depend 
on the differentiable structure  of $M$. In any case,
there is no obvious way to read off these invariants from
homotopy invariants of the $M$.

Nonetheless, homotopy invariants do impose some 
important relations between these invariants. 
One such relation arises from the 
intersection form 
\begin{eqnarray*}
\smile  :
H^{2}(M, {\mathbb R})\times H^{2}(M, {\mathbb R})	
 & \longrightarrow & ~~~~ {\mathbb R}  \\
	( ~ [\phi ] ~ , ~ [\psi ] ~) ~~~~~ & 
	\mapsto  & \int_{M}\phi \wedge \psi 
\end{eqnarray*}
which may be  diagonalized as
	    $$\left[ 
	  \begin{array}{rl}  \underbrace{
	 \begin{array}{ccc}
	   		 1 &  &	  \\
	   		  &	\ddots &   \\
	   		  &	 & 1
	   	  \end{array}}_{b_{+}(M)}
	   	   & 
	   		   \\ 
	  
	 {\scriptstyle b_{-}(M)}  \!
	    \left\{\begin{array}{r}
	    \\
	    \\
	    \\
	    \\
	    \end{array}
	    \right. \! \! \! \! \! \! \! \! \! \! \! \! \! \! \! 
	    &\begin{array}{ccc}
	   		 -1	&  &   \\
	   		  &	\ddots &   \\
	   		  &	 & -1
	   	  \end{array}
	    \end{array} 
	    \right] $$
by choosing a suitable basis for the de Rham cohomology 
$H^{2}(M, {\mathbb R})$. 
The numbers
$b_{\pm} (M)$ are independent of the
choice of basis, and so are  oriented 
homotopy invariants of $M$. Their difference
$$\tau (M) = b_{+}(M) - b_{-}(M),$$
 is called the {\em signature} of $M$. 
The Hirzebruch signature theorem \cite{hirz,st}
asserts that this invariant is expressible as
a curvature integral: 
\begin{equation}
	\tau (M)= \frac{1}{12\pi^2}\int_M \left(|W^+|^2
- |W^-|^2  \right) d\mu .
	\label{sig}
\end{equation}
Here the curvatures, norms $|\cdot |$, and
volume form $d\mu$ are,  of course, those of the 
Riemannian metric $g$, but the entire point is that
the answer is independent of {\em which} metric we use. 
Thus  $(\inf \|W_{+}\|)^{2}=(\inf \|W_{-}\|)^{2}+ 12\pi^{2}\tau (M)$.

Another relationship between the curvature L$^{2}$-norms under
consideration is given \cite{aw,st} by the $4$-dimensional case of 
the generalized Gauss-Bonnet theorem; this 
asserts that the Euler characteristic 
$$\chi (M) = 2-2b_{1}(M) + b_{2}(M)$$
is given by 
\begin{equation}
	\chi (M)= \frac{1}{8\pi^2}\int_M \left(|W^+|^2
+ |W^-|^2 + \frac{s^2}{24} -\frac{|\stackrel{\circ}{r}|^2}{2}
\right)d\mu .
	\label{eul}
\end{equation}
This then gives rise to inequalities concerning  
$\inf  \|W_{\pm}\|$, $\inf  \|s\|$, and $\inf  \| \ro \|$
which are imposed by  the homotopy type of $M$.

Not long ago,  the problem of calculating 
invariants such as $\inf \| s\|$ would have simply  been
considered intractable. However, there has been 
a remarkable amount of recent progress on these issues.
The key event in this regard 
was   Witten's introduction \cite{witten}  of 
 the so-called Seiberg-Witten invariants,
which  display an unexpected  relationship between 
 Donaldson's polynomial invariants \cite{don} and Riemannian
geometry on $M$. 
In particular, $\inf \| s\|$ turns out 
{\em not} to be a   homeomorphism invariant, but is  
nonetheless exactly calculable 
 for a huge class 
of 4-manifolds, including all  complex 
algebraic surfaces  \cite{lno,lky}.
In the next section, we will give a brief introduction
to Seiberg-Witten theory, and then show how
it allows one to estimate certain linear combinations
of $\|s\|$ and $\|W_{+}\|$.

\section{Seiberg-Witten Theory}
\label{weylest}

Let $(M,g)$ be a  compact oriented Riemannian $4$-manifold.
On any contractible open subset  $U\subset M$, one can define 
Hermitian vector  bundles 
$$\begin{array}{rl}
	{\mathbb C}^{2}\rightarrow & {\mathbb S}_{\pm}|_{U}  \\
	 & \downarrow  \\
	 & U\subset M
\end{array}$$
called spin bundles,  characterized by the 
fact that their determinant  line bundles $\wedge^{2}{\mathbb S}_{\pm}$
are canonically trivial and that their projectivizations
$$\begin{array}{rc}
	\bcp_{1}\rightarrow & {\mathbb P}({\mathbb S}_{\pm} ) \\
	 & \downarrow  \\
	 & M
\end{array}$$
are exactly the unit 2-sphere bundles  $S(\Lambda^{\pm})$.
As one passes between open subset $U$ and 
$U'$, however, the corresponding locally-defined spin bundles 
are not quite canonically isomorphic; instead, 
there are two equally  `canonical'
isomorphisms, differing by 
 a  sign. 
Because of this, one cannot generally define the bundles
${\mathbb S}_{\pm}$ globally on $M$; manifolds on which 
this can be done are called {\em spin}, and are characterized by 
the vanishing of 
the Stiefel-Whitney class $w_{2}=w_{2}(TM)\in H^{2}(M,\ZZ_{2})$. 
However, 
one can always find Hermitian complex line bundles $L\to M$
 with first Chern class $c_{1}=c_{1}(L)$ satisfying  
 \begin{equation}
	c_{1}\equiv 
w_{2} ~~\bmod 2.
	\label{redux}
\end{equation}
 Given such a line bundle, one can then construct
 Hermitian vector bundles ${\mathbb V}_{\pm}$ with 
$${\mathbb P}({\mathbb V}_{\pm})= S(\Lambda^{\pm})$$
by formally setting 
$${\mathbb V}_{\pm}= {\mathbb S}_{\pm}\otimes L^{1/2},$$
 because the sign problems encountered in consistently 
 defining the 
 transition functions of ${\mathbb S}_{\pm}$ are exactly 
 canceled by those associated with trying to find 
 consistent square-roots of the transition functions of 
 $L$. 
 
 The isomorphism class $\mathfrak c$ of  such a 
 choice of ${\mathbb V}_{\pm}$ is called 
  a {\em spin$^{c}$ structure}
on $M$.
 The cohomology group $H^{2}(M,{\mathbb Z})$ acts freely and
transitively on the spin$^{c}$ structures by tensoring ${\mathbb V}_{\pm}$
with complex line bundles. 
  Each spin$^{c}$
 structure has a first Chern  class
$c_{1}:= c_{1}(L)= c_{1}({\mathbb V}_{\pm}) \in H^{2}(M,{\mathbb Z})$ 
 satisfying (\ref{redux}),
 and the  
$H^{2}(M,{\mathbb Z})$-action 
on spin$^{c}$ structures  
induces the action 
$$c_{1}\mapsto c_{1}+2\alpha,$$ $\alpha \in H^{2}(M,{\mathbb Z})$,
on first Chern classes. Thus, if 
$H^{2}(M,{\mathbb Z})$  has trivial 2-torsion  --- as can always 
be arranged by replacing $M$ with a finite cover ---
the spin$^{c}$ structures are precisely in one-to-one correspondence 
with the set of cohomology classes $c_{1} \in H^{2}(M,{\mathbb Z})$ 
satisfying (\ref{redux}). 

To make this discussion more concrete, suppose that
$M$ admits an almost-complex structure.
 Any given almost-complex
structure can be deformed to an almost complex structures $J$ 
 which is compatible with $g$
in the sense that $J^*g=g$. Choose such a $J$, and consider
the  rank-2 complex vector bundles 
\begin{eqnarray} {\mathbb V}_+&=& \Lambda^{0,0}\oplus \Lambda^{0,2}\label{spl}\\
{\mathbb V}_-&=&\Lambda^{0,1}.\nonumber \end{eqnarray}
These are precisely the twisted spinor bundles of the 
spin$^{c}$ structure obtained by taking 
$L$ to be the anti-canonical line bundle
$\Lambda^{0,2}$ of the almost-complex structure. 
A spin$^{c}$ structure $\mathfrak c$  arising  in this way will
be said to be of {\em almost-complex type}.
These are exactly the  spin$^{c}$ structures 
for which 
$$c_{1}^{2}= (2\chi + 3\tau )(M).$$

On a spin manifold, the spin bundles ${\mathbb S}_{\pm}$
carry natural connections induced by the Levi-Civita connection
of the given Riemannian metric $g$. 
On a spin$^{c}$ manifold, however, there is not a 
natural unique choice
of  connections on ${\mathbb V}_{\pm}$. Nonetheless,  since we formally have 
${\mathbb V}_{\pm}= {\mathbb S}_{\pm}\otimes L^{1/2}$, every 
 Hermitian connection $A$ on  $L$ induces associated
  Hermitian connections 
$\nabla_{A}$ on ${\mathbb V}_{\pm}$. 

On the other hand, 
there is a canonical isomorphism 
$\Lambda^{1}\otimes {\mathbb C}=\mbox{Hom }({\mathbb 
S}_{+}, {\mathbb S}_{-})$, so that 
$\Lambda^{1}\otimes {\mathbb C} \cong \mbox{Hom }({\mathbb V}_{+}, 
{\mathbb V}_{-})$ for any spin$^{c}$ structure,
and this induces a canonical homomorphism 
$$\cdot ~: \Lambda^{1}\otimes  {\mathbb V}_{+}\to {\mathbb V}_{-}$$
called {\em Clifford multiplication}. Composing these
operations allows us to define a so-called {\em twisted  Dirac 
operator} 
$$D_{A}:\Gamma ({\mathbb V}_{+})\longrightarrow \Gamma ({\mathbb V}_{-})$$
by $D_{A}\Phi = \nabla_{A}\cdot \Phi$.

 %

For any spin$^{c}$ structure, we have already noted that 
there is a canonical diffeomorphism 
${\mathbb P}({\mathbb V}_{+})\stackrel{\simeq}{\to} S(\Lambda^{+})$. 
In polar coordinates, we now use this 
to define  the angular part of a unique continuous map 
$$\sigma : {\mathbb V}_{+} \to \Lambda^{+} $$
with 
$$|\sigma (\Phi ) | = \frac{1}{2\sqrt{2}}|\Phi |^{2}.$$
This map is actually real-quadratic on each fiber of ${\mathbb V}_{+}$;
indeed, assuming our spin$^{c}$ 
structure is induced by a complex structure $J$, then, in terms of 
(\ref{spl}), $\sigma$ is explicitly given by  
$$\sigma (f, \phi)=(|f|^2-|\phi|^2) \frac{\omega}{4}+ \Im m 
(\bar{f}\phi),$$
where $f\in \Lambda^{0,0}$, $\phi \in \Lambda^{0,2}$, and where
$\omega (\cdot, \cdot ) = g(J\cdot, \cdot )$ is the 
associated  2-form of $(M,g,J)$.  On a deeper level,
$\sigma$ directly  arises from the fact that 
${\mathbb V}_{+}={\Bbb S}_{+}\otimes L^{1/2}$, while 
$\Lambda^{+}\otimes {\mathbb C}= \odot^{2}{\Bbb S}_{+}$.
For this reason, $\sigma$ is 
 is invariant under parallel transport.
  
We are now in a position to introduce the
 Seiberg-Witten equations 
\begin{eqnarray} D_{A}\Phi &=&0\label{drc}\\
 F_{A}^+&=&i \sigma(\Phi) ,\label{sd}\end{eqnarray}
where the  unknowns are a Hermitian connection $A$ on $L$
and a  section $\Phi$ of ${\mathbb V}_+$.
Here 
$F_{A}^+$  is the self-dual part of
the curvature of $A$,
and so is a
 purely imaginary 2-form.

For many $4$-manifolds, it turns out that  
there is a solution of the Seiberg-Witten equations 
for each metric. 
Let us  introduce some convenient 
 terminology \cite{K} to describe this situation.

\begin{defn}
Let $M$ be a smooth compact oriented $4$-manifold
with $b_{+}\geq 2$, and suppose that $M$
 carries a   spin$^{c}$ structure $\mathfrak c$
for which the  Seiberg-Witten 
equations (\ref{drc}--\ref{sd})
have a solution for every Riemannian  metric $g$ on $M$. 
Then  the first Chern class $c_{1}\in H^{2}(M,\ZZ )$ 
of $\mathfrak c$ will be called a   
{\em monopole class}. 
\end{defn}

This definition is useful in practice primarily because there
are topological arguments which lead to   
 the existence of solutions the Seiberg-Witten equations.
For example \cite{witten}, if $\mathfrak c$  is  a spin$^c$ structure 
of almost-complex type, 
 then   the  {\em Seiberg-Witten invariant}   
${\mathcal S \mathcal W}_{\mathfrak c}(M)$ 
can be defined as the number of solutions,
modulo gauge transformations and  counted with orientations, of a  
generic perturbation 
\begin{eqnarray*} D_{A}\Phi &=&0  \\
 iF^+_A+\sigma (\Phi ) &=& \phi  \end{eqnarray*}
 of  (\ref{drc}--\ref{sd}), 
where $\phi$ is a smooth self-dual 2-form. 
If $b_{+}(M)\geq 2$, this integer is 
independent of the metric $g$; and if it is non-zero,
the first Chern class $c_{1}$ of $\mathfrak c$ is then
a monopole class. Similar things are also true when $b_{+}(M)=1$, 
although the story becomes  rather more complicated.  

Now, via the Hodge theorem, every 
 Riemannian metric $g$ on	$M$	 determines
  a direct	sum	decomposition
  $$H^{2}(M, {\Bbb	R})	=
  {\mathcal  H}^{+}_{g}\oplus {\cal	H}^{-}_{g},$$
  where ${\cal	H}^{+}_{g}$	 (respectively,	${\mathcal  H}^{-}_{g}$)
 	consists of	those cohomology
  classes for which the   harmonic	 representative	 is	
 	self-dual  (respectively,  anti-self-dual).
 	Because the restriction of the intersection form to
 	${\mathcal  H}^{+}_{g}$ (respectively, ${\cal	H}^{-}_{g}$)
 	is positive (respectively, negative) definite, and because
 	these subspaces are mutually orthogonal with respect to
 	the intersection pairing, the dimensions of these spaces
 	are exactly the invariants
 	$b_{\pm}$ defined  in \S \ref{prim}. If 
 	the first Chern class $c_{1}$ of the spin$^{c}$ structure
 	$\mathfrak c$ is now decomposed as 
 	$$c_{1}= c_{1}^{+}+c_{1}^{-},$$
 	where $c_{1}^{\pm}\in {\mathcal  H}^{\pm}_{g}$, 
 	we get  the important inequality	
  \begin{equation}
 	   \int_{M}|\Phi |^{4}d\mu \geq	32\pi^{2}(c_{1}^{+})^{2}
 	   \label{harm}
  \end{equation}
  because (\ref{sd}) tells us that 
  $2\pi c_{1}^{+}$	is the harmonic	part of	
  $-\sigma	(\Phi )$.

 Many of the most remarkable consequences of  Seiberg-Witten 
 theory stem from the fact that the equations (\ref{drc}--\ref{sd})  
imply the Weitzenb\"ock formula
\begin{equation}\label{wb}
0= 4\nabla^*\nabla \Phi +  s \Phi +|\Phi|^2\Phi ,
\end{equation}
where $s$ denotes the scalar curvature of $g$, and where
we have introduced the abbreviation $\nabla_{A}=\nabla$.
Taking the inner product with $\Phi$, it follows that
\begin{equation}
	0=2\Delta |\Phi |^{2} +4|\nabla \Phi |^{2} + s |\Phi |^{2}+ |\Phi 
|^{4} .
	\label{zero}
\end{equation}
If we  multiply (\ref{zero}) by 
 $|\Phi |^{2}$ and integrate,  we  have
$$0 = \int_{M}\left[ 2\left| ~d|\Phi |^{2} ~\right|^{2} + 4 |\Phi |^{2}|\nabla 
\Phi |^{2}+ s|\Phi |^{4}+ |\Phi |^{6}
\right] d\mu_{g} ,$$
so that 
\begin{equation}
	\int (-s) |\Phi |^{4} d\mu \geq 4 \int |\Phi |^{2} |\nabla \Phi 
|^{2} d\mu  + \int  |\Phi |^{6} d\mu .
	\label{two}
\end{equation}
This leads  \cite{lric} to the following
 curvature estimate:

\begin{thm} \label{L2}
Let $M$ be a   smooth compact  
oriented  4-manifold with monopole class $c_{1}$. 
Then  every  
Riemannian metric $g$ on $M$ 
satisfies 
\begin{equation}
	\int_{M}\left({\frac{2}{3}s-2\sqrt{\frac{2}{3}}|W_{+}|}\right)^{2} 
	d\mu 
\geq 32\pi^{2}(c_{1}^{+})^{2}, 
	\label{central}
\end{equation}
where $c_1^+$  is the self-dual part of $c_{1}$ with respect to 
$g$.
\end{thm}
\begin{proof} 
The first step is to  prove the inequality 
 \begin{equation}
 	V^{1/3}\left( 
 	\int_{M}\left|\frac{2}{3}s_{g}-2\sqrt{\frac{2}{3}}|W_{+}|
 	\right|^{3}d\mu 
 \right)^{2/3}\geq 32\pi^{2} (c_{1}^{+})^{2},
 	\label{crux}
 \end{equation}
 where $V=\mv (M,g)= \int_{M}d\mu_{g}$ is the total volume of $(M,g)$. 

 Any self-dual  2-form $\psi$ on any oriented 4-manifold satisfies 
 the 
Weitzenb\"ock formula \cite{bourg}
$$(d+d^{*})^{2}\psi = \nabla^{*}\nabla \psi - 2W_{+}(\psi , 
\cdot ) + \frac{s}{3} \psi.$$
It follows that 
$$
\int_{M}(-2W_{+})(\psi , \psi ) 	d\mu \geq 
\int_{M}(-\frac{s}{3})|\psi |^{2}~d\mu -
 \int_{M} |\nabla \psi |^{2}
 ~d\mu . 	
$$
However, 
$$ |W_{+}|_{g}|\psi |^{2}\geq - \sqrt{\frac{3}{2}}W_{+}(\psi , \psi )$$
simply because $W_{+}$ is trace-free.
Thus
 $$
  \int_{M}2\sqrt{\frac{2}{3}}|W_{+}| |\psi |^{2} 	d\mu \geq 
\int_{M}(-\frac{s}{3})|\psi |^{2}~d\mu -
\int_{M} |\nabla \psi |^{2}
 ~d\mu , 
  $$ 
  and hence
   $$
  -\int_{M}(\frac{2}{3}s-2\sqrt{\frac{2}{3}}|W_{+}|)|\psi |^{2} 	d\mu \geq 
\int_{M}(-s)|\psi |^{2}~d\mu -
\int_{M} |\nabla \psi |^{2}
 ~d\mu . 
  $$ 
On the other hand, the  particular self-dual 2-form $\varphi = \sigma (\Phi 
)=-iF_{A}^{+}$ satisfies 
\begin{eqnarray*}
	|\varphi |^{2} & = & \frac{1}{8}|\Phi |^{4}  ,\\
	|\nabla \varphi |^{2} & \leq  & \frac{1}{2} |\Phi |^{2}|\nabla \Phi |^{2} .
\end{eqnarray*}
 Setting $\psi = \varphi$,  
  we thus have
  $$
   -\int_{M}(\frac{2}{3}s-2\sqrt{\frac{2}{3}}|W_{+}|)|\Phi |^{4} 	d\mu \geq 
\int_{M}(- s)|\Phi |^{4}~d\mu -
 4\int_{M} |\Phi |^{2} |\nabla \Phi |^{2}
 ~d\mu .
  $$
 But   (\ref{two})
 tells us that 
 $$\int_{M}(- s)|\Phi |^{4}~d\mu -
 4\int_{M} |\Phi |^{2} |\nabla \Phi |^{2}
 ~d\mu \geq \int_{M}|\Phi |^{6}~d\mu ,$$
 so we obtain 
 \begin{equation}
 	 -\int_{M}(\frac{2}{3}s-2\sqrt{\frac{2}{3}}|W_{+}|)|\Phi |^{4} d\mu \geq  
 \int_{M}|\Phi |^{6}~d\mu .
 	\label{voici}
 \end{equation}
By the H\"older inequality, we thus have
$$\left( \int \left|\frac{2}{3}s-2\sqrt{\frac{2}{3}}|W_{+}|
\right|^{3}d\mu \right)^{1/3}
\left( \int |\Phi |^{6} d\mu \right)^{2/3}
\geq  
 \int|\Phi |^{6}~d\mu ,$$
 Since  the H\"older inequality also tells us 
 that
 $$
 \int|\Phi |^{6}~d\mu  \geq V^{-1/2}\left(\int |\Phi |^{4}d\mu \right)^{3/2} ,
 $$
 we thus have 
$$ V^{1/3}\left( \int_{M}\left|\frac{2}{3}s-2\sqrt{\frac{2}{3}}|W_{+}|
\right|^{3}d\mu 
 \right)^{2/3}\geq 
 \int |\Phi |^{4}d\mu
 \geq 32\pi^{2} (c_{1}^{+})^{2},$$
where the last inequality is exactly (\ref{harm}).  
 This completes the first part of the proof. 
 
 Next, we  observe that any smooth conformal $\gamma$   class on 
 any oriented $4$-manifold  contains a $C^{2}$ metric such that 
 $s-\sqrt{6}|W_{+}|$ is  constant. Indeed, as observed 
 by Gursky \cite{G1},
 this readily follows from the standard 
 proof of the Yamabe problem. The main point is that 
   the curvature expression
 $${\mathfrak S}_{g}= s_{g}-\sqrt{6}|W_{+}|_{g}$$
 transforms under conformal changes $g\mapsto \hat{g}= u^{2}g$
 by the rule 
 $${\mathfrak S}_{\hat{g}}=u^{-3}\left(6\Delta_{g} + {\mathfrak 
 S}_{g}\right)u ,$$
 just like the ordinary scalar curvature $s$. 
  We will actually use this only 
 in the negative case, where the proof is technically the simplest,
 and simply repeats\footnote{However, since  $|W_{+}|$ is generally only 
 Lipschitz continuous,  the minimizer generally only has regularity 
 $C^{2,\alpha}$ in the vicinity of a zero of $W_{+}$.} 
 the arguments of Trudinger \cite{trud}.
  
 The conformal class $\gamma$ of a given metric $g$
 thus always contains 
a metric $g_{\gamma}$ for which 
$\frac{2}{3}s-2\sqrt{\frac{2}{3}}|W_{+}|$ is constant.
But since the existence of solutions of the Seiberg-Witten
equations precludes the possibility that we might have 
$s_{g_{\gamma}}>0$, this  constant is necessarily non-positive. 
We thus have 
$$\int_{M}\left(\frac{2}{3}s_{g_{\gamma}}-2\sqrt{\frac{2}{3}}|W_{+}|_{g_{\gamma}}
 \right)^{2}d\mu_{g_{\gamma}} =
V^{1/3}_{g_{\gamma}}\left( \int_{M}\left|(\frac{2}{3}s_{g_{\gamma}}-
{2}\sqrt{\frac{2}{3}}|W_{+}|_{g_{\gamma}} \right|^{3}d\mu_{g_{\gamma}} 
 \right)^{2/3},$$
 so that
 $$
 \int_{M}\left(\frac{2}{3}s_{g_{\gamma}}-2\sqrt{\frac{2}{3}}|W_{+}|_{g_{\gamma}}
 \right)^{2}d\mu_{g_{\gamma}} \geq 
  32\pi^{2} (c_{1}^{+})^{2}. 
 $$
 Thus we at least have the desired  $L^{2}$ 
 estimate for a specific metric $g_{\gamma}$
 which is conformally related to the given metric $g$.
 
Let us now compare the left-hand side with analogous 
expression for 
 the given metric $g$. To do so, 
we express $g$ in the form $g=u^{2}g_{\gamma}$,
where $u$ is a positive $C^{2}$ function, and observe that 
\begin{eqnarray*}
	\int_{M} \left(\frac{2}{3}s_{g}-2\sqrt{\frac{2}{3}}|W_{+}|_{g}
 \right) u^{2}d\mu_{g_{\gamma}}  & = & \frac{2}{3}\int 
 {\mathfrak S}_{g}u^{2}d\mu_{g_{\gamma}}  \\
	 & = & \frac{2}{3}\int u^{-3}\left(6\Delta_{g_{\gamma}}u+
	 {\mathfrak S}_{g_{\gamma}}u\right) u^{2}d\mu_{g_{\gamma}}
	 \\&=& \frac{2}{3}\int \left(-6 u^{-2} |du|^{2}_{g_{\gamma}}+
	 {\mathfrak S}_{g_{\gamma}}\right)  d\mu_{g_{\gamma}}
	 \\ &\leq & \frac{2}{3}\int 
	 {\mathfrak S}_{g_{\gamma}}  d\mu_{g_{\gamma}}\\
	 &=& \int_{M} \left(\frac{2}{3}s_{g_{\gamma}}
	 -2\sqrt{\frac{2}{3}}|W_{+}|_{g_{\gamma}}\right) 
	 d\mu_{g_{\gamma}} .
\end{eqnarray*}
Applying Cauchy-Schwarz, we thus have
\begin{eqnarray*}
 -V^{1/2}_{g_{\gamma}}\left(\int \left(\frac{2}{3}s_{g}
 -2\sqrt{\frac{2}{3}}|W_{+}|_{g}
 \right)^{2} d\mu_{g}\right)^{1/2}	
 & \leq  & 
 \int_{M} \left(\frac{2}{3}s_{g}-2\sqrt{\frac{2}{3}}|W_{+}|_{g}
 \right) u^{2}d\mu_{g_{\gamma}}
 \\&\leq &\int_{M} \left(\frac{2}{3}s_{g_{\gamma}}
 -2\sqrt{\frac{2}{3}}|W_{+}|_{g_{\gamma}}\right) 
	 d\mu_{g_{\gamma}}  \\
	 & = &  -
	 V^{1/2}_{g_{\gamma}}\left(\int 
	 \left(\frac{2}{3}s_{g_{\gamma}}-2\sqrt{\frac{2}{3}}|W_{+}|_{g_{\gamma}}
	  \right)^{2}d\mu_{g_{\gamma}}\right)^{1/2} , 
\end{eqnarray*}
and hence  
$$\int_{M}\left(\frac{2}{3}s_{g}-2\sqrt{\frac{2}{3}}|W_{+}|_{g} 
\right)^{2}d\mu_{g} \geq 
\int_{M}\left(\frac{2}{3}s_{g_{\gamma}}-2\sqrt{\frac{2}{3}}|W_{+}|_{g_{\gamma}}
 \right)^{2}d\mu_{g_{\gamma}} ,
$$
exactly as claimed. 
  \end{proof}

Notice that we can   rewrite the inequality (\ref{central})  as
 $$\left\|\frac{2}{3}s-2\sqrt{\frac{2}{3}}|W_{+}|\right\|\geq 4\sqrt{2}\pi 
 |c_{1}^{+}|,$$
 where $\|\cdot\|$ denotes the $L^{2}$ norm with respect to $g$.
 Dividing by $\sqrt{24}$ and applying the triangle inequality, 
 we thus have 

\begin{cor} \label{encore}
Let $M$ be a   smooth compact  
oriented  4-manifold with monopole class $c_{1}$.
Then  every  
Riemannian metric $g$ on $M$ 
satisfies  
 \begin{equation}
 	\frac{2}{3}\|\frac{s}{\sqrt{24}}\| + {\frac{1}{3}}\| W_{+}\|\geq 
 	\frac{2\pi}{\sqrt{3}} 
 |c_{1}^{+}|. 
 	\label{trngl}
 \end{equation}
\end{cor}

Inequality (\ref{central}) actually belongs to a family of 
related estimates:    

\begin{thm} \label{fin}
Let $M$ be a   smooth compact  
oriented  4-manifold with monopole class $c_{1}$, 
and let $\delta\in [0,\frac{1}{3}]$ be a constant.  
Then  every  
Riemannian metric $g$ on $M$ 
satisfies 
\begin{equation}
	\int_{M}\left[(1-\delta) s- \delta \sqrt{24}|W_{+}|\right]^{2} 
	d\mu 
\geq 32\pi^{2}(c_{1}^{+})^{2}, 
	\label{gen}
\end{equation}
\end{thm}
\begin{proof}
Inequality (\ref{two}) implies
\begin{equation}
	\int (-s) |\Phi |^{4} d\mu \geq  \int  |\Phi |^{6} d\mu .
	\label{deux}
\end{equation}
On the other hand, inequality (\ref{voici}) 
asserts that 
$$	 -\int_{M}(\frac{2}{3}s-2\sqrt{\frac{2}{3}}|W_{+}|)|\Phi |^{4} d\mu \geq  
 \int_{M}|\Phi |^{6}~d\mu .
$$
Now multiply (\ref{deux}) by $1-3\delta$, multiply (\ref{voici})
by $3\delta$, and add. The result is 
\begin{equation}
	\int \left[(1-\delta) s-\delta \sqrt{24}|W_{+}|\right]
	 |\Phi |^{4} d\mu \geq  \int  |\Phi |^{6} d\mu .
	\label{voila}
\end{equation} 
Applying the same H\"older inequalities as before, we now obtain 
$$ V^{1/3}\left( \int_{M}\left| (1-\delta) s-\delta \sqrt{24}|W_{+}|
\right|^{3}d\mu 
 \right)^{2/3}\geq 
 \int |\Phi |^{4}d\mu
 \geq 32\pi^{2} (c_{1}^{+})^{2}.$$
 
 Passage from this $L^{3}$ estimate to the desired 
 $L^{2}$ estimate is then accomplished by the same
 means as before: every conformal class contains a metric
 for which $(1-\delta) s-\delta \sqrt{24}|W_{+}|$
 is constant, and this metric minimizes 
 	$$\int_{M} \left[(1-\delta) s-\delta \sqrt{24}|W_{+}| \right]^{2} 
	d\mu $$
	among metrics in its  conformal class.  
\end{proof}

Rewriting  (\ref{gen})  as
 $$\left\|(1-\delta) s-\delta \sqrt{24}|W_{+}|\right\|\geq 4\sqrt{2}\pi 
 |c_{1}^{+}|,$$
 dividing by $\sqrt{24}$, 
 and
 applying the triangle inequality, we thus have 

\begin{cor} \label{bis} 
Let $M$ be a   smooth compact  
oriented  4-manifold with monopole class $c_{1}$.
Then  every  
Riemannian metric $g$ on $M$ 
satisfies  
 \begin{equation}
 	(1-\delta)\|\frac{s}{\sqrt{24}}\| + \delta \| W_{+}\|
 	\geq \frac{2\pi}{\sqrt{3}} 
 |c_{1}^{+}| 
 	\label{gentl}
 \end{equation}
 for every $\delta \in [0,\frac{1}{3}]$. 
\end{cor}

The  $\delta = 0$ version of (\ref{gen}) is implicit in the work of 
Witten \cite{witten}; it was later made explicit in  
\cite{lpm}, where it was also shown  that
equality holds for $\delta =0$ iff $g$ is a K\"ahler metric
of constant, non-positive scalar curvature. But indeed, since  
$\sqrt{24}|W_{+}|\equiv |s|$ for any K\"ahler manifold of real
dimension $4$,   metrics of this kind saturate  (\ref{gen})
for each value of $\delta$. Conversely:

\begin{prop} Let $\delta\in [0,\frac{1}{3})$ be a fixed constant.
If $g$ is a metric such that 
equality holds in (\ref{gen}), then
 $g$ is K\"ahler, and has constant
scalar curvature.
\end{prop}
\begin{proof}
Equality in (\ref{gen})  implies  equality  in 
(\ref{voila}). However,  $(1-3\delta )$ times  inequality (\ref{two}) 
plus $3\delta$ times inequality (\ref{voici}) 
reads 
$$
	\int \left[(1-\delta) s-\delta \sqrt{24}|W_{+}|\right]
	 |\Phi |^{4} d\mu \geq  \int  |\Phi |^{6} d\mu +
	4 (1-3\delta ) \int |\Phi |^{2}|\nabla \Phi |^{2}d\mu .
$$
Equality in (\ref{gen}) therefore implies that 
$$0 =\frac{1}{2} \int |\Phi |^{2}|\nabla \Phi |^{2}d\mu
\geq 
\int |\nabla \varphi |^{2}d\mu ,$$
forcing  the $2$-form $\varphi$ to be parallel. 
If $\varphi \not\equiv 0$, we conclude that the metric is 
K\"ahler, and the constancy of $s$ then follows from 
the Yamabe portion of the argument. 

On the other hand, since $b_{+}(M)\geq 2$ and  $c_{1}$ is a monopole class, 
$M$ does not admit any metrics of positive scalar curvature.
If  $\varphi\equiv 0$ and (\ref{gen}) is saturated, 
one can therefore show that  $(M,g)$ is 
$K3$ or $T^{4}$ with a  Ricci-flat K\"ahler metric. 
The details are left as an exercise for the
interested reader. 
\end{proof}

When $\delta = \frac{1}{3}$, the above argument breaks down. 
However, a metric $g$ can  saturate (\ref{central})  only
if equality holds in (\ref{harm}), and this forces the
self-dual $2$-form $\varphi= \sigma (\Phi )$ to be 
harmonic. Moreover, the relevant H\"older inequalities
would  also have to be saturated, forcing $\varphi$ to have
constant length. This forces $g$ to be 
 {\em almost-K\"ahler}, in the sense that 
there is an orientation-compatible orthogonal 
almost-complex structure for which the associated
$2$-form is closed. For details, see \cite{lric}. 

It is reasonable to ask whether the inequalities  (\ref{gen}) 
and (\ref{gentl}) continue to hold when $\delta > 1/3$.
This issue will be addressed in \S \ref{sharp}.

\section{Einstein Metrics}

Recall that a smooth Riemannian metric $g$ is said to be {\em  Einstein} 
if its Ricci curvature $r$ is a constant multiple of 
the metric: 
$$r=\lambda g .$$
Not every   4-manifold admits such metrics.
A necessary condition for the existence of an Einstein metric
on a compact oriented 4-manifold is that the 
Hitchin-Thorpe inequality $2\chi  (M) \geq  3 |\tau (M)|$
must hold \cite{tho,hit,bes}. Indeed, (\ref{sig}) and (\ref{eul}) tell
us that 
$$		(2\chi \pm 3\tau ) (M) 
=  \frac{1}{4\pi^2}\int_M \left(
  \frac{s^2}{24}+ 2|W_\pm|^2-\frac{|\stackrel{\circ}{r}|^2}{2}
\right)d\mu.
$$
The Hitchin-Thorpe inequality follows, since
the integrand is non-negative when $\ro =0$.  
This argument, however, treats the 
scalar and Weyl contributions as `junk' terms, 
about which one knows nothing except that they are
non-negative. We now remedy this by 
invoking the estimates of \S \ref{weylest}. 

\begin{prop}
Let $M$ be a   smooth compact  
oriented  4-manifold with monopole class $c_{1}$. 
Then every  metric $g$ on $M$ satisfies
$$
\frac{1}{4\pi^{2}}\int_{M}\left( \frac{s_{g}^{2}}{24} + 
2|W_{+}|_{g}^{2}\right) d\mu_{g} \geq \frac{2}{3} 
(c_{1}^{+})^{2} .
$$
If $c_{1}^{+}\neq 0$, moreover, 
equality can only hold if $g$ is
almost-K\"ahler, with almost-K\"ahler class proportional 
to $c_{1}^{+}$. 
\end{prop}
\begin{proof}
 We begin  begin with inequality (\ref{trngl})  
 $$
\frac{2}{3}\|\frac{s}{\sqrt{24}}\| + {\frac{1}{3}}\| W_{+}\|\geq 
 	\frac{2\pi}{\sqrt{3}} 
 |c_{1}^{+}|,$$
 and  elect  to interpret the left-hand side as the dot product
 $$(\frac{2}{3} , {\frac{1}{3\sqrt{2}}}) \cdot \left(\|\frac{s}{\sqrt{24}}\|, 
 \sqrt{2}\| W_{+}\| \right) $$
 in ${\mathbb R}^{2}$. Applying  Cauchy-Schwarz, we thus have 
 $$ \left((\frac{2}{3})^{2} + ({\frac{1}{3\sqrt{2}}})^{2}\right)^{1/2}
 \left(\int_{M}(\frac{s^{2}}{24}+2|W_{+}|^{2})d\mu \right)^{1/2}
 \geq \frac{2}{3}\|\frac{s}{\sqrt{24}}\| + {\frac{1}{3}}\| W_{+}\| .$$
 Thus
 $$\frac{1}{2}\int_{M}(\frac{s^{2}}{24}+2|W_{+}|^{2})d\mu \geq 
 \frac{4\pi^{2}}{3}(c_{1}^{+})^{2},$$
and hence   
 $$\frac{1}{4\pi^{2}}\int_{M}\left( \frac{s_{g}^{2}}{24} + 
2|W_{+}|_{g}^{2}\right) d\mu_{g} \geq \frac{2}{3} 
(c_{1}^{+})^{2} ,$$
as claimed. 

In the equality case, $\varphi$ would be a closed
self-dual form of  constant norm, so $g$ would be 
almost-K\"ahler unless $\varphi\equiv 0$. 
\end{proof}

To give some  concrete  applications, we now  
focus on the case of complex surfaces. 

\begin{prop}\label{west}
Let $(X,J_{X})$ be a compact complex surface with 
$b_{+}>1$, and let $(M,J_{X})$ be the complex surface
obtained from $X$ by blowing up $k>0$ points. 
Then any Riemannian metric $g$ on the  $4$-manifold   
$$M= X\# k \overline{\bcp}_{2}$$
satisfies
$$
\frac{1}{4\pi^{2}}\int_{M}\left( \frac{s_{g}^{2}}{24} + 
2|W_{+}|_{g}^{2}\right) d\mu_{g} > \frac{2}{3} 
(2\chi + 3\tau)(X) .
$$
\end{prop}
\begin{proof}
Let
 $c_{1}(X)$ denote the first Chern class of 
 the given complex structure $J_{X}$, and,  by a standard  abuse of 
 notation, let $c_{1}(X)$ also denote the  
 pull-back class of this class to $M$.
 If 
  $E_{1}, \ldots , E_{k}$ are the Poincar\'e duals of
  the exceptional divisors in $M$ introduced by blowing up,
  the  complex structure $J_{M}$ has Chern class 
  $$c_{1}(M)= c_{1}(X)-\sum_{j=1}^{k}E_{j}.$$
 By  a result of  Witten \cite{witten}, this 
  is a monopole class of $M$. 
  However, there are self-diffeomorphisms of $M$
  which act on $H^{2}(M)$ in a manner such that 
  \begin{eqnarray*}
  	c_{1}(X) & \mapsto  &  c_{1}(X) \\
  	E_{j} & \mapsto & \pm E_{j}
  \end{eqnarray*}
  for any choice of signs we like. 
 Thus $$c_{1}= c_{1}(X) + \sum_{j=1}^{k} (\pm  E_{j})$$
  is a monopole class on $M$ for each choice of signs. 
  We now fix  our choice of signs so that
  $$[c_{1}(X)]^{+}\cdot (\pm E_{j})\geq 0,$$
  for each $j$, with respect to the decomposition induced by 
  the given metric $g$.  
 We then have 
 \begin{eqnarray*}
 	(c_{1}^{+})^{2} & = & \left([c_{1}(X)]^{+} + 
 	\sum_{j=1}^{k} (\pm	E_{j}^{+})\right)^{2}  \\
 	 & = & ([c_{1}(X)]^{+})^{2}+ 2 \sum_{j=1}^{k} [c_{1}(X)]^{+}\cdot
 	 (\pm  E_{j})
 + (\sum_{j=1}^{k} (\pm E_{j}^{+}))^{2}  \\
 	 & \geq  & ([c_{1}(X)]^{+})^{2}  \\
 	 & \geq  & (2\chi + 3\tau ) (X).
  \end{eqnarray*}
  This shows that 
  $$
\frac{1}{4\pi^{2}}\int_{M}\left( \frac{s_{g}^{2}}{24} + 
2|W_{+}|_{g}^{2}\right) d\mu_{g} \geq \frac{2}{3} 
(2\chi + 3\tau)(X) .
$$

If equality held, $g$ would be almost-K\"ahler, with 
 almost-K\"ahler class 
$[\omega ]$  proportional to $c_{1}^{+}$. 
On the other hand, we would also have 
 $[c_{1}(X)]^{+}\cdot  E_{j}=0$,
 so it  would then follow that   
  $[\omega ]\cdot E_{j}=0$ for all $j$. However, the
 Seiberg-Witten invariant would  be  non-trivial for a spin$^{c}$ 
 structure with $c_{1}(\tilde{L})= c_{1}(L)-2 (\pm E_{1})$,
and a celebrated theorem of Taubes \cite{taubes3} would then force 
the homology class $E_{j}$ to be represented by a pseudo-holomorphic
$2$-sphere in the symplectic manifold $(M,\omega )$. But the 
(positive!) area of this
sphere with respect to  $g$ would 
then be  exactly $[\omega ]\cdot E_{j}$, contradicting the 
observation 
that $[\omega ]\cdot E_{j}=0$. 
\end{proof}

\begin{thm}\label{blowup}
Let $(X,J_{X})$ be a compact complex surface with $b_{+}>1$, and 
let $(M,J_{M})$ be obtained from $X$ by blowing up 
$k$  points. Then the smooth compact $4$-manifold $M$ 
does not admit any Einstein metrics if 
 $k \geq 
\frac{1}{3}c_{1}^{2}(X)$. 
\end{thm}
\begin{proof}
We may assume that $(2\chi + 3\tau ) (X) > 0$, since otherwise
the result follows from the Hitchin-Thorpe inequality. 

Now 
$$(2\chi + 3\tau )(M)=\frac{1}{4\pi^{2}}\int_{M}\left( \frac{s_{g}^{2}}{24} + 
2|W_{+}|_{g}^{2} -\frac{|\stackrel{\circ}{r}|^{2}}{2}\right) d\mu_{g}$$
for any metric on $g$ on $M$.  If $g$ is an Einstein metric,
the trace-free part $\stackrel{\circ}{r}$ of the Ricci curvature
vanishes, and  we then have
\begin{eqnarray*}
(2\chi + 3\tau )(X) -k 	 & = & (2\chi + 3\tau )(M)  \\
	 & = &  \frac{1}{4\pi^{2}}\int_{M}\left( \frac{s_{g}^{2}}{24} + 
2|W_{+}|_{g}^{2}\right) d\mu_{g} \\
	 & > & \frac{2}{3} (2\chi + 3\tau )(X)  \\
\end{eqnarray*}
 by Proposition \ref{west}.  
 If $M$ carries an Einstein metric, it therefore follows 
 that 
 $$ \frac{1}{3}(2\chi + 3\tau )(X) > k .$$
 The claim thus follows by contraposition. 
\end{proof}

\begin{xpl} Let $X\subset {\bcp}_{4}$ 
be the intersection of two cubic hypersurfaces in general 
position. Since the canonical class on $X$ is exactly the
hyperplane class, $c_{1}^{2}(X)= 1^{2}\cdot 3 \cdot 3= 9$.
 Theorem
\ref{blowup}
therefore tells us that if we  
blow up $X$  at $3$ points, the  resulting $4$-manifold
$$M = X\# 3 \overline{\bcp}_{2}$$
does not admit Einstein metrics.

But now consider the {\em  Horikawa surface} $N$ 
obtained as a ramified double  cover of 
the blown-up projective plane 
${\mathbb C \mathbb P}_{2}\# \overline{\mathbb C \mathbb P}_{2}$ 
 branched over the (smooth) proper transform $\hat{C}$ of the singular 
 curve $C$ given by 
 $$x^{10}+y^{10}+ z^{6}(x^{4}+y^{4})=0$$
 in  the complex projective plane, where the singular point
 $[0:0:1]$ of $C$ is the point at which we blow up ${\mathbb C \mathbb P}_{2}$.
  By the  Freedman classification of $4$-manifolds \cite{freedman},
  both of these complex surfaces  
 are homeomorphic to 
 $$11{\mathbb C \mathbb P}_{2}\#53\overline{\mathbb C \mathbb P}_{2}.$$
However, $N$ has $c_{1}<0$, and 
so admits a K\"ahler-Einstein metric by the Aubin/Yau theorem
\cite{aubin,yau}. 
Thus, although $M$ and $N$ are homeomorphic, one admits
Einstein metrics, while the other doesn't.
\end{xpl}

\begin{xpl} Let $X\subset {\bcp}_{3}$ be a hypersurface of 
degree $6$. Since the canonical class on $X$ is twice the
hyperplane class, $c_{1}^{2}(X)= 2^{2}\cdot 6= 24$.  Theorem
\ref{blowup}
therefore tells us that if we  
blow up $X$  at $8$ points, the  resulting $4$-manifold
$$M = X\# 8 \overline{\bcp}_{2}$$
does not admit Einstein metrics. 

However, the Freedman classification can be used to show that  
$M$ is homeomorphic to the   Horikawa surface $N$ 
obtained as a ramified double  cover of 
${\mathbb C \mathbb P}_{1}\times {\mathbb C \mathbb P}_{1}$ 
 branched at a generic curve of bidegree 
 $(6,12)$; indeed, both of these  complex surfaces
 are 
 homeomorphic to 
 $$21{\mathbb C \mathbb P}_{2}\#93\overline{\mathbb C \mathbb P}_{2}.$$
However, this $N$ also
 admits a K\"ahler-Einstein metric, even though the existence of
 Einstein metric is obstructed on $M$. 
\end{xpl}
 %
 %
 %
 %

 	\begin{xpl}	
 	 Let $X\subset {\bcp}_{3}$ be a	hypersurface of	
  degree $10$.	Since the canonical	class on $X$ is	six	times the
  hyperplane class, $c_{1}^{2}(X)=	6^{2}\cdot 10= 360$.  Theorem
  \ref{blowup}
  therefore tells us that if we  
  blow	up $X$	at $120$ or	more points, the  resulting	$4$-manifold
  does	not	admit Einstein metrics.	In particular, 
  this	assertion applies to 
  $$M = X\# 144 \overline{\bcp}_{2}.$$

  Now let $N$ be obtained from	 
 	${\mathbb C	\mathbb	P}_{1}\times {\mathbb C	\mathbb	P}_{1}$	
 	as a ramified double cover 
 	branched at	a generic curve	of bidegree	
 	$(8,58)$. Both $M$ and $N$ are then	simply connected, and have 
 	$c_1^{2}= 216$ and $p_{g}=84$; and both	
 	are	therefore homeomorphic to 
 	$$129 {\mathbb C \mathbb P}_{2}\# 633
 	\overline{\mathbb C	\mathbb	   P}_{2}.$$
 	But	again, $N$ has $c_{1}<0$, and so 
 	admits a K\"ahler-Einstein metric, even	though
 	$M$	does not admit an Einstein metric of any kind whatsoever.
 	
 	In most	respects, this example is much like	the	
 	previous examples. However,	
 	this choice	of $N$ is not a	Horikawa surface, but instead 
 	sits well away from	the	Noether	line \cite{bpv}	of 
 	complex-surface	geography.
 	\end{xpl}

Infinitely many such examples  can be constructed using the above
techniques, and the interested reader might wish to explore their
geography. 

It should be noted  that  
Theorem \ref{blowup} is the direct descendant  of 
an analogous result in \cite{lno}, where
scalar curvature estimates alone were used to 
obtain an obstruction when $k \geq \frac{2}{3} c_{1}^{2}(X)$. It 
was later pointed out by Kotschick \cite{kot} that this 
suffices to imply the existence of homeomorphic pairs
consisting of an Einstein manifold and a 
$4$-manifold which does not admit
 Einstein metrics.   
 An intermediate step between \cite{lno} and 
 Theorem \ref{blowup} may be found in \cite{lebweyl},
 where  cruder Seiberg-Witten estimates of Weyl curvature
  were  used to obtain 
  an obstruction for $k \geq \frac{25}{57} c_{1}^{2}(X)$.

  \section{How Sharp are the Estimates?}
  \label{sharp}

 The estimates we have described in \S \ref{weylest} 
 are optimal in the sense that equality is achieved  
 for K\"ahler metrics of constant negative scalar curvature,
 with the standard orientation and spin$^{c}$ structure.
 In this section, we will attempt to probe the 
 limits of these estimates by considering metrics 
 of precisely this type, but with {\em non-standard}
 choices of orientation and spin$^{c}$ structure. 
 
 One interesting class of $4$-manifolds  which admit 
 constant-scalar-curvature K\"ahler metrics are 
 the complex surfaces with ample canonical line bundle.
 In terms of complex-surface classification \cite{bpv}, these are
 precisely those minimal surfaces of general type which 
 do not contain  $\bcp_{1}$'s of self-intersection 
 $-2$.
 The ampleness of the canonical line bundle
 is often written as $c_{1}<0$, meaning that 
 $-c_{1}$ is a K\"ahler class. 
 A celebrated result of Aubin/Yau \cite{aubin,yau} 
 guarantees that there is a unique K\"ahler-Einstein
 metric on $M$, compatible with the given 
 complex structure, and with K\"ahler
 class   
 $[\omega] = -c_{1}= H^{1,1}(M, \RR )$.
 The scalar curvature of such a metric is, of course,  
 a negative constant; indeed,  $s=-\dim_{\RR}M=-4$.

 Now if $M$ is a compact complex manifold
 without holomorphic vector fields,
 the set of K\"ahler classes which are representable 
 by metrics of constant scalar curvature is open
 \cite{fusch,ls} 
 in $H^{1,1}(M, \RR )$. On the other hand, 
 a manifold with $c_{1}< 0$ never carries a 
 non-zero holomorphic vector field, so it follows
 that a complex surface with ample canonical line bundle 
 will carry lots of constant-scalar-curvature  K\"ahler metrics 
 which are non-Einstein if $b_{-}=h^{1,1}-1$ is non-zero. 
 However, one might actually hope to find such metrics
 even in those K\"ahler classes which are far from the
 anti-canonical class. This expectation may be codified as follows:

\begin{conj}\label{hip} 
Let $M$ be any compact complex surface with $c_{1}  < 0$.
Then every K\"ahler class $[\omega]\in H^{1,1}(M,\RR )$ contains a
unique K\"ahler metric of constant  scalar curvature.
\end{conj}

The uniqueness clause was recently proved by 
X.-X. Chen \cite{xxchen}, using ideas due to 
Donaldson  and Semmes.
A direct continuity-method attack on conjecture has also been
explored by S.-R. Simanca.

Let us now narrow our discussion to a very special class of
complex surfaces. 

 \begin{defn}
A {\em Kodaira fibration} is  a holomorphic
submersion $\varpi: M\to B$  from a compact complex surface to a 
compact complex curve, such that the  base $B$ and fiber 
$F_{z}= \varpi^{-1}(z)$ both have genus $\geq 2$. 
If $M$ admits
such a fibration $\varpi$, we will say that is a 
{\em Kodaira-fibered surface}. 
\end{defn}

The underlying $4$-manifold $M$ of a
Kodaira-fibered surface is 
a fiber bundle over $B$, with fiber $F$. We thus have \cite{span} 
a long exact sequence
$$\cdots \to \pi_{k}(F) \to \pi_{k} (M) \to \pi_{k}(B) \to 
\pi_{k-1}(F)\to \cdots$$
of homotopy groups, so that 
$M$ is a $K(\pi , 1)$. Thus, any $2$-sphere
in $M$ is homologically trivial, and so has 
self-intersection $0$; in particular,  the complex surface $M$ cannot contain
any $\bcp_{1}$'s of self-intersection $-1$ or 
 $-2$. On the other hand, $M$ is of general type, 
 so the above implies that  $c_{1}(M) < 0$.
 Kodaira-fibered surfaces thus provide us with 
 an interesting testing-ground for Conjecture \ref{hip}. 

Now the product $B\times F$ of two complex curves of
genus $\geq 2$ is certainly  Kodaira fibered, but such  
a product also admits orientation-reversing diffeomorphisms, and 
so has  signature
$\tau=0$. However,
as was first observed by 
Kodaira \cite{kodf}, one can  construct examples 
with $\tau > 0$ by taking {\em branched covers} of 
products; cf.  \cite{atkodf,bpv}. 
 For example, let $B$ be a 
 curve of genus $3$ with a holomorphic 
 involution $\iota : B\to B$ without fixed points; 
 one may visualize such an involution as a $180^{\circ}$ rotation 
 of a $3$-holed doughnut about an axis which passes though
 the middle hole, without meeting the doughnut. 
 Let $f: C\to B$ be the unique $64$-fold unbranched cover
 with $f_{*}[\pi_{1}(C)]= \ker [\pi_{1}(B) \to H_{1}(B, \ZZ_{2})]$;
 thus $C$ is a complex curve of genus $129$. Let 
 $\Sigma \subset C\times B$ be the union of 
 the graphs of $f$ and $\iota\circ f$. Then
 the homology class of $\Sigma$ is divisible by $2$. We may
 therefore construct a ramified double cover $M\to B\times C$
 branched over $\Sigma$. The projection $M\to B$
 is then a Kodaira fibration, with fiber $F$ of genus $321$.
 The projection $M\to C$ is also a Kodaira fibration, 
 with fiber of genus $6$. The signature of this example is 
  $\tau (M) =256$, and so  coincidentally equals one-tenth of 
  its Euler characteristic
  $\chi (M) = 2560$.

Now, more generally, let $M$ be  any Kodaira-fibered surface
with $\tau > 0$, and let  
$\varpi : M\to B$ be a Kodaira fibration.
Let $p$ denote the the genus of $B$, and let
$q$ denote the genus of a fiber $F$ of $\varpi$.
Indulging in a standard notational abuse, let us also 
use $F$ to denote the 
 Poincar\'e dual
 of the homology class of the fiber. Since 
  $F$ can be represented in de Rham cohomology 
  by the pull-back of an area
form on $B$, this $(1,1)$-class is positive semi-definite. On the other hand, 
$-c_{1}$ is a K\"ahler class on $M$, and so it follows that 
$$
[\omega_{\epsilon}] = 2(p-1) F- \epsilon c_{1}
$$ 
is a K\"ahler class on $M$ for any 
$\epsilon > 0$. If Conjecture \ref{hip}
is true, there must therefore exist a 
K\"ahler metric $g_{\epsilon}$ 
on $M$ of constant scalar curvature
with K\"ahler class $[\omega_{\epsilon}]$.
Let us explore the global geometric invariants of this 
putative metric.

The metric in question, being K\"ahler,  would have total scalar curvature 
$$\int s_{g_{\epsilon}} d\mu_{g_{\epsilon}} = 4\pi c_{1}\cdot 
[\omega_{\epsilon}]= -4\pi (\chi +\epsilon c_{1}^{2}) (M)$$
and total volume
$$\int d\mu_{g_{\epsilon}} = \frac{[\omega_{\epsilon}]^{2}}{2}=
\frac{\epsilon}{2}(2\chi +\epsilon c_{1}^{2}) (M).$$
The assumption that $s_{g_{\epsilon}}=\mbox{const}$
would thus imply  that
\begin{eqnarray*}\|s\|^{2}=
\int s^{2}_{g_{\epsilon}}d\mu_{g_{\epsilon}}	 & = & 
\frac{32\pi^{2}}{\epsilon} 
\frac{(\chi +\epsilon c_{1}^{2})^{2}}{2\chi +\epsilon c_{1}^{2}}  \\
	 & = & 16 \pi^{2} \frac{\chi}{\epsilon} \left[ {1} + 
	 (3 + \frac{9}{2}\varrho ) {\epsilon}+ 
O(\epsilon^{2} )\right],
\end{eqnarray*}
where we have set
$$\varrho = \frac{\tau (M)}{\chi (M)}  ~~.$$
Since a  K\"ahler metric on a complex surface satisfies 
$|W_{+}|^{2}\equiv s^{2}/24$,  
we would also consequently 
have 
\begin{eqnarray*}
	\int |W_{+}|^{2}_{g_{\epsilon}}d\mu_{g_{\epsilon}} & = &
	 \frac{1}{24} \int s^{2}_{g_{\epsilon}}d\mu_{g_{\epsilon}}  \\
	 & = & \frac{2}{3} \pi^{2}\frac{\chi}{\epsilon}  \left[ {1} 
	 + (3 + \frac{9}{2}\varrho ) {\epsilon}+ 
O(\epsilon^{2} )\right].
\end{eqnarray*}
It would thus follow that 
\begin{eqnarray*}
\|W_{-}\|^{2}=	\int |W_{-}|^{2}_{g_{\epsilon}}d\mu_{g_{\epsilon}} & = &
	 -12\pi^{2}\tau (M) 
	 +\int |W_{+}|^{2}_{g_{\epsilon}}d\mu_{g_{\epsilon}}  \\
	 & = & \frac{2}{3} \pi^{2} \frac{\chi}{\epsilon} 
	 \left[{1} + (3 - \frac{27}{2}\varrho ) {\epsilon}+ 
O(\epsilon^{2} )\right].
\end{eqnarray*}

On the other hand, there are symplectic forms on $M$ which 
are compatible with the {\em non-standard} orientation 
of $M$; for example, the cohomology class $F+\epsilon c_{1}$
is represented by such forms if $\epsilon$ is sufficiently small. 
A celebrated theorem of Taubes \cite{taubes} therefore 
tells us that the reverse-oriented version
$\overline{M}$ of $M$ has a non-trivial Seiberg-Witten 
invariant  \cite{leung,jpik,kot2}. The relevant 
spin$^{c}$ structure on $\overline{M}$ is of 
almost-complex type, and its first Chern class, which
we will denote by $\bar{c}_{1}$, is given by
$$\bar{c}_{1}= c_{1}+4(p-1)F.$$
Of course, the conjugate almost-complex structure, 
with first Chern class $-\bar{c}_{1}$,is  also a
monopole class of  $\overline{M}$, and $\overline{M}$
will have yet other monopole  classes if, for example, 
$M$ admits more than one Kodaira fibration and $\tau (M) \neq 0$. 

Now recall that (\ref{gentl}) asserts that
$$
(1-\delta) \|\frac{s}{\sqrt{24}}\| + 
\delta \| W_{+}\|\geq \frac{2\pi}{\sqrt{3}}
 |c_{1}^{+}|
$$
for all $\delta \in [0,\frac{1}{3}]$. 
One would like to know  whether this inequality 
 might also hold, quite generally, for some value of 
$\delta > \frac{1}{3}$. 
In order to find out, we  apply this inequality to $\overline{M}$ with the
above monopole class.  Rewriting the  inequality 
with respect to the {\em complex} orientation of
$M$, we then get 
\begin{equation}
(1-\delta) \|\frac{s}{\sqrt{24}}\| + 
\delta \| W_{-}\|\geq \frac{2\pi}{\sqrt{3}}
 |{\bar{c}_{1}}^{-}| ,
	\label{foobar}
\end{equation}
and it is this inequality we shall now use to probe the
limits of the theory.

Relative to any K\"ahler metric with K\"ahler class
$[\omega_{\epsilon}]$, one has
\begin{eqnarray*}
\bar{c}_{1}^{+}	 & = & 
\frac{\bar{c}_{1}\cdot [\omega_{\epsilon}]}{[\omega_{\epsilon}]^{2}}
[\omega_{\epsilon}]  \\
	 & = & \frac{[c_{1}+4(p-1)F]\cdot
	  [2(p-1)F-\epsilon c_{1}]}{[\omega_{\epsilon}]^{2}}
[\omega_{\epsilon}] \\
	 & = & -\frac{(\chi + 3\epsilon \tau)}{[\omega_{\epsilon}]^{2}}
	 [\omega_{\epsilon}] ,
\end{eqnarray*}
so that 
\begin{eqnarray*}
|\bar{c}_{1}^{+}|^{2}	 & = & 
\frac{(\chi + 3\epsilon \tau)^{2}}{[\omega_{\epsilon}]^{2}}  \\
	 & = &\frac{1}{\epsilon}
	  \frac{(\chi + 3\epsilon \tau)^{2}}{2\chi +\epsilon c_{1}^{2}}  \\
	 & = & \frac{\chi }{2\epsilon}\left[
	 1- (1-\frac{9}{2}\varrho)\epsilon + O(\epsilon^{2} ) 
	 \right] .
\end{eqnarray*}
Now since $\bar{c}_{1}$ is an almost-complex structure on 
$\overline{M}$, we have 
$$|{\bar{c}_{1}}^{-}|^{2}-|\bar{c}_{1}^{+}|^{2}= 2\chi - 3\tau, $$
so 
\begin{eqnarray*}
|{\bar{c}_{1}}^{-}|^{2}	 & = & (2\chi - 3\tau) + \frac{\chi }{2\epsilon}\left[
	 1- (1-\frac{9}{2}\varrho)\epsilon + O(\epsilon^{2} ) 
	 \right]
 \\
	 & = & \frac{\chi}{2}(4 - 6 \varrho) + \frac{\chi }{2\epsilon}\left[
	 1- (1-\frac{9}{2}\varrho)\epsilon + O(\epsilon^{2} ) 
	 \right]  \\
	 & = &   \frac{\chi }{2\epsilon}\left[
	 1+ (3-\frac{3}{2}\varrho)\epsilon + O(\epsilon^{2} ) 
	 \right].
\end{eqnarray*}

After dividing by $\pi \sqrt{2\chi/3\epsilon}$, the inequality (\ref{foobar}) 
would  thus read 
$$
(1-\delta) \sqrt{ 
{1} + (3 + \frac{9}{2}\varrho ){\epsilon}+ 
O(\epsilon^{2})} + 
\delta \sqrt{
	{1} + (3 - \frac{27}{2}\varrho ){\epsilon}+ 
O(\epsilon^{2} )}\geq \sqrt{
	 1+ (3-\frac{3}{2}\varrho)\epsilon + O(\epsilon^{2} ) 
	 } .
$$
Dropping the terms of order $\epsilon^{2}$, we would thus have 
$$
(1-\delta) \left[
{1} + (\frac{3}{2} + \frac{9}{4}\varrho ){\epsilon}\right] + 
\delta \left[
	{1} + (\frac{3}{2} - \frac{27}{4}\varrho ){\epsilon}\right]
	\geq 
	 1+ (\frac{3}{2}-\frac{3}{4}\varrho)\epsilon ,	 
$$
so that, upon collecting terms, we would obtain 
$$3\varrho \epsilon \geq 9 \varrho \epsilon \delta.$$
Taking $\varrho = \tau/\chi$ to be positive, and 
noting that $\epsilon$  is positive by construction,
this shows that  Conjecture \ref{hip} would imply that 
$$\frac{1}{3} \geq \delta ,$$
or in other words that (\ref{trngl}) is optimal.
We have thus proved the following result: 

\begin{thm}
Either 
\begin{itemize}
\item the estimate (\ref{trngl}) is optimal; or else 
\item Conjecture \ref{hip} is false. 
\end{itemize}
\end{thm}

 %

\end{document}